\documentclass[11pt]{amsart}

\usepackage{amsfonts,epsfig}
\usepackage{latexsym}
\usepackage{amssymb}
\usepackage{amsmath}
\usepackage{amsthm}
\usepackage{graphics}
\usepackage[all]{xy}
\usepackage[T2A]{fontenc}
\usepackage{multirow}
\usepackage{array, booktabs, ctable}
\usepackage{enumerate}
\usepackage{booktabs}
\usepackage{longtable}
\usepackage{multirow}
\usepackage[pagebackref]{hyperref}
\usepackage{xcolor}

\newtheorem{theorem}{Theorem}

\newcommand{\cC}{{\mathcal{C}}}

\newcommand{\Q}{{\mathbb{Q}}}
\newcommand{\Qbar}{\overline{\mathbb Q}}
\newcommand{\Z}{{\mathbb{Z}}}

\DeclareMathOperator{\cm}{CM} 
\DeclareMathOperator{\tors}{tors} 

\newcommand{\fcm}{{\mathfrak{cm}}}

\title[Torsion of rational elliptic curves with CM over cubic fields]{Torsion growth over cubic fields\\ of rational elliptic curves with complex multiplication}

\author{Enrique Gonz\'alez--Jim\'enez}
\address{Departamento de Matem{\'a}ticas, Universidad Aut{\'o}noma de Madrid, Madrid, Spain}
\email{enrique.gonzalez.jimenez@uam.es}
\urladdr{http://matematicas.uam.es/~enrique.gonzalez.jimenez}

\thanks{The author was partially  supported by the grant PGC2018-095392-B-I00.}

\begin{document}
\date{July 16, 2019}

\maketitle

\begin{abstract} 
This article is a contribution to the project of classifying the torsion growth of elliptic curve upon base-change. In this article we treat the case of elliptic curve defined over the rationals with complex multiplication. For this particular case, we give  a description of the possible torsion growth over cubic fields and a completely explicit description of this growth in terms of some invariants attached to a given elliptic curve.
\end{abstract}

\section{Introduction}
The arithmetic of elliptic curves is one of the most fascinating areas in Number Theory or Algebraic Geometry. Let $E$ be an elliptic curve defined over a number field $K$, then the Mordell-Weil Theorem asserts that the set of $K$-rational points on $E$, denoted by $E(K)$, forms a finitely generated abelian group. The subgroup of points of finite order, denoted by $E(K)_{\tors}$, is called the torsion subgroup and it is well known that is isomorphic to $\cC_n\times\cC_m$ for some positive integers $n,m$, where $\cC_n=\Z/n\Z$ denotes the cyclic group of order $n$. The study of torsion subgroups is had been treated for several active people last years. Thanks to Merel \cite{Merel}, it is known that given a positive integer $d$, the set $\Phi(d)$  of possible groups (up to isomorphism) that can appear as the torsion subgroup $E(K)_{\tors}$, where $K$ runs through all number fields $K$ of degree $d$ and $E$ runs through all elliptic curves over $K$, is finite. Only the cases $d=1$ and $d=2$ are known (by \cite{Mazur1978}; and \cite{Kamienny92,KenkuMomose88} respectively).\\
This paper focuses on a particular approach concerning torsion growth: we are interested in studying how does the torsion subgroup of an elliptic curve defined over $\Q$ change when we consider the elliptic curve over a number field of degree $d$. Note that if $E$ is an elliptic curve defined over $\Q$ and $K$ a number field such that the torsion of $E$  grows from $\Q$ to $K$, then of course the torsion of $E$ also grows from $\Q$ to any extension of $K$. We say that the torsion growth over $K$ is primitive if $E(K')_{\tors}\subsetneq E(K)_{\tors}$ for any subfield $K'\subsetneq K$.\\
We introduce some useful definition for the sequel:
\begin{itemize}
\item Let $\Phi_\Q(d)$ be the set of possible groups (up to isomorphisms) that can appear as the torsion subgroup over a number field of degree $d$, of an elliptic curve defined over $\Q$. 
\item Fixed $G \in \Phi(1)$, let $\Phi_\Q(d,G)$ be the subset of $\Phi_\Q(d)$ such that $E$ runs through all elliptic curves over $\Q$ such that $E(\Q)_{\tors}\simeq G$.
\item  Let $E$ be an elliptic curve defined over $\Q$ and d a positive integer. We denote by $\mathcal{H}_{\Q}(d,E)$ the multiset of groups $H$ such that there exist $K/\Q$, an extension of degree dividing $d$, with $H\simeq E(K)_{\tors}\ne E(\Q)_{\tors}$ and the torsion growth in $K$ is primitive. Note that we are allowing the possibility of two (or more) of the torsion subgroups $H$ being isomorphic if the corresponding number fields $K$ are not isomorphic. We let $\mathcal{H}_{\Q}(d)$ denote the set of $\mathcal{H}_{\Q}(d,E)$ as $E$ runs over all  elliptic curves defined over $\Q$. Finally, for any $G\in \Phi(1)$ we define $\mathcal{H}_{\Q}(d,G)$ as the set of multisets $\mathcal{H}_{\Q}(d,E)$ where $E$ runs over all the elliptic curve defined over $\Q$ such that $E(\Q)_{\tors}\simeq G$. Denote by $h_\Q(d)$ the maximum of the cardinality of $S$ when $S\in \mathcal{H}_{\Q}(d)$, then $h_\Q(d)$ gives the maximum number of field extension of degrees dividing $d$ where there is primitive torsion growth. 
\end{itemize}
The sets $\Phi_\Q(d)$, $\Phi_\Q(d,G)$ and $\mathcal{H}_{\Q}(d,G)$, for any $G \in \Phi(1)$, have been completely classified for $d=2,3,5,7$ and for any positive integer $d$ whose prime divisors are greater than $7$ (cf. \cite{Najman16,GJT14,GJT15,GJNT16, GJ17,GJN16}). The set $\Phi_\Q(4)$ is known \cite{Chou16,GJN16}. The other sets have been treated for  $d=4$ in \cite{GJLR18} and $d=6$ in \cite{HGJ18}. 
We denote by  $\Phi^{\cm}(d)$, $\Phi^{\cm}_\Q(d)$,\\ $\Phi^{\cm}_\Q(d,G)$, $\mathcal{H}^{\cm}_{\Q}(d,G)$, which are defined as the analogues above sets but restricting to elliptic curves with complex multiplication (CM).

The set $\Phi^{\cm}(1)$ was determined by Olson \cite{Olson74}:
$$
\Phi^{\cm}(1)=\left\{ \cC_1\,,\,  \cC_2\,,\,  \cC_3\,,\,  \cC_4\,,\,  \cC_6\,,\, \cC_2\times\cC_2\right\}.
$$
The quadratic cases by M\"uller  et al \cite{MSZ89} and the cubic case by several authors headed by Zimmer \cite{FSWZ90,PWZ97}:
$$
\Phi^{\cm}(3) = \Phi^{\cm}(1)\cup  \left\{\,\cC_9\,,\,\cC_{14}\,\right\}.
$$
Recently, Clark et al. \cite{Clark2014} have computed the sets $\Phi^{\cm}(d)$, for $4\le d\le 12$.
\

Restricting to elliptic curve with complex multiplication defined over $\Q$ we obtain the following results:

\begin{theorem}\label{teo1} $\Phi_{\Q}^{\cm}(3)=\Phi^{\cm}\left( 3\right)$.
\end{theorem}

\begin{theorem}\label{teo2}
Let be $G \in \Phi^{\cm}(1)$. Then 
\begin{itemize}
\item If $G\in \{\cC_4,\cC_6,\cC_2\times\cC_2\}$ then $\Phi^{\cm}_\Q \left(3,G \right)=\{G\}$. In particular, $\mathcal{H}^{\cm}_{\Q}(3,G)=\emptyset$.
\item If $G\in \{\cC_1,\cC_2,\cC_3\}$ then the sets $\Phi^{\cm}_\Q \left(3,G \right)$ and $\mathcal{H}^{\cm}_{\Q}(3,G)$ are the following:\\[2mm]
\begin{center}
\begin{tabular}{|c|c|l|}
\hline
$G$ & $ \Phi^{\cm}_\Q \left(3,G \right)\setminus\{G\}$ & $\mathcal{H}^{\cm}_{\Q}(3,G)$\\
\hline\hline
\multirow{3}{*}{$\cC_1$}  & \multirow{3}{*}{$\left\{\, {\cC_2,\cC_3,\cC_6}\, \right\}$} & {$\cC_2$}     \\
\cline{3-3}
& & {$\cC_6$}   \\
\cline{3-3}
& & {$\cC_2\,,\cC_3$}   \\
\hline
\multirow{2}{*}{$\cC_2$}  & \multirow{2}{*}{$\left\{\, {\cC_6}\,,\, \cC_{14}\, \right\}$} & {$\cC_6$}     \\
\cline{3-3}
& & {$ \cC_{14}$}    \\
\hline
\multirow{2}{*}{$\cC_3$}  & \multirow{2}{*}{$\left\{\, \cC_6\,,\, \cC_{9}\, \right\}$} & {$\cC_6$}     \\
\cline{3-3}
& & ${\cC_6}\,,\, \cC_{9}$   \\
\hline
\end{tabular}
\end{center}
\end{itemize}
In particular, $h^{\cm}_\Q(3)=2$.
\end{theorem}

Our aim in this paper is to go further. More precisely, one we have given a description of the possible torsion growth over cubic fields we are going to give a completely explicit description of this growth in terms of some invariants attached to a given elliptic curve. The case of quadratic growth is solved in \cite{GJCM2}. In an ongoing paper \cite{GJCMn} we will solve the problem for number fields of low degree. 
\begin{theorem}\label{teo3}
Table \ref{d3} gives an explicit description of torsion growth over cubic fields of any elliptic curve defined over $\Q$ with CM depending only in its corresponding CM-invariants (see \S\ref{appendix} for the definition).
\end{theorem}

\begin{table}[ht!]
 \caption{Explicit description of torsion growth over cubic fields of elliptic curves defined over $\Q$ with complex multiplication}\label{d3}
\begin{tabular}{|c|c|c|l|c|}
\hline
 $\mathfrak{cm}$ & $k$ such that $E=E^k_{\mathfrak{cm}}$ & $G\simeq  E(\Q)_{\tors}$ & $\mathcal{H}_\Q(E,3)$ &  cubics $\Q(\alpha)$\\
\hline\hline
 \multirow{9}{*}{$3$} & $1$ & $\cC_6$ & $-$  & $-$ \\
\cline{2-5}
 & $16$&   \multirow{3}{*}{$\cC_3$}  & {$\cC_6,\,\cC_9$}  &  $\sqrt[3]{2},\,\, \alpha^3-3\alpha-1=0$ \\
\cline{2-2}\cline{4-5}
 & $-432$&    & {$\cC_6$}  &  $\sqrt[3]{2}$ \\
\cline{2-2}\cline{4-5}
 & $r^2 \,\, (r\ne\pm 1,\pm 4)$ &   & {$\cC_6$}  & $\sqrt[3]{k}$  \\
\cline{2-5}
 &$-27$  &  \multirow{2}{*}{$\cC_2$}  & {$\cC_6$}  &  $\sqrt[3]{2}$ \\
\cline{2-2}\cline{4-5}
 &{$r^3\,\,(r\ne 1,-3)$}  &  & $-$  & $-$  \\
\cline{2-5}
 & $-108$&   \multirow{3}{*}{$\cC_1$}  & {$\cC_6$}  &  $\sqrt[3]{2}$ \\
\cline{2-2}\cline{4-5}
 & $-3 r^2\,\, (r\ne \pm 6)$ & & $\cC_2,\,\cC_3$ & $\sqrt[3]{3r^2},\,\,\sqrt[3]{12r^2}$\\
\cline{2-2}\cline{4-5}
 & $\ne r^2,r^3,-3r^2$ & & $\cC_2$ & $\sqrt[3]{k}$\\
\cline{1-5}
  \multirow{3}{*}{$12$} & $1$ & $\cC_6$ & $-$  & $-$ \\
\cline{2-5}
  & $-3$ & \multirow{2}{*}{$\cC_2$} & $\cC_6$  & $\sqrt[3]{2}$\\
\cline{2-2}\cline{4-5}
  & $\ne 1,-3$ & & $-$  & $-$\\
\cline{1-5}
  \multirow{3}{*}{$27$} & $1$ & $\cC_3$ & {$\cC_6,\,\cC_9$}  &  $\sqrt[3]{2},\,\, \alpha^3-3\alpha-1=0$ \\
\cline{2-5}
  & $-3$ & \multirow{2}{*}{$\cC_1$} & $\cC_2,\,\cC_3$  & $\sqrt[3]{2},\,\sqrt[3]{3}$\\
\cline{2-2}\cline{4-5}
  & $\ne 1,-3$ & & $\cC_2$  & $\sqrt[3]{2}$\\
\hline
 \multirow{3}{*}{$4$} & $4$ &{$\cC_4$} & $-$  &$-$ \\
\cline{2-5}
 & $-r^2$ & {$\cC_2\times\cC_2$} & $-$  & $-$\\
\cline{2-5}
 & $\ne 4,-r^2$ & {$\cC_2$} & $-$  & $-$\\
\cline{1-5}
  \multirow{2}{*}{$16$} & $1,2$ &{$\cC_4$} & $-$  & $-$\\
 \cline{2-5}
  &$\ne 1,2$ & {$\cC_2$} & $-$  &$-$ \\
\hline
   \multirow{2}{*}{$7$} & $-7$ &  \multirow{2}{*}{$\cC_{2}$} & $\cC_{14}$  & $\alpha^3+\alpha^2-2\alpha-1=0$\\
\cline{2-2}\cline{4-5}
 & $\ne -7$ &   & $-$  & $-$\\
\cline{1-5}
  \multirow{2}{*}{$28$} & $7$ &  \multirow{2}{*}{$\cC_{2}$} & $\cC_{14}$  & $\alpha^3+\alpha^2-2\alpha-1=0$\\
\cline{2-2}\cline{4-5}
 & $\ne 7$ &  & $-$  & $-$\\
 \hline
 $8$ &$-$ & {$\cC_2$} & $-$  & $-$ \\
 \hline
$11$ &$-$ & $\cC_1$ & $\cC_2$  &  $\alpha^3-\alpha^2+\alpha+1=0$ \\
 \hline
$19$ &$-$ & $\cC_1$ & $\cC_2$  & $\alpha^3-\alpha^2+3\alpha-1=0$ \\
 \hline
 $43$ & $-$& $\cC_1$ & $\cC_2$  &$\alpha^3-\alpha^2-\alpha+3=0$ \\
 \hline
 $67$ & $-$& $\cC_1$ & $\cC_2$  & $\alpha^3-\alpha^2-3\alpha+5=0$\\
 \hline
 $163$ & $-$ & $\cC_1$ & $\cC_2$  & $\alpha^3-8\alpha-10=0$\\
 \hline
 \end{tabular}

\end{table}

{\em Notation:}  Given an elliptic curve $E:y^2=x^3+Ax+B$, $A,B\in K$, and a number field $K$, we denote by $j(E)$ its $j$-invariant, by $\Delta(E)$ the discriminant of that short Weierstrass model, and by $E(K)_{\tors}$ the torsion subgroup of the Mordell-Weil group of $E$ over $K$. For a positive integer $n$, we denote by $\cC_n=\Z/n\Z$ the cyclic group of order $n$.

\section{Proof of the Theorems}

\subsection{Preliminaries}  Let $E$ be an elliptic curve and $n$ a positive integer. Denote by $E[n]$ the set of points on $E$ of order dividing $n$. The $x$-coordinates of the points on $E[n]$ correspond to the roots of the $n$-division polynomial $\Psi_n(x)$ of $E$ (cf. \cite[\S 3.2]{Washington}). By abuse of notation, in this paper we use $\Psi_n(x)$ to denote the primitive $n$-division polynomial of $E$, that is, the classical $n$-division polynomial factors by the $m$-division polynomials of $E$ for proper factors $m$ of $n$. Then $\Psi_n(x)$ is characterized by the property that its roots are the $x$-coordinates of the points of exact order $n$ of $E$. In particular if $E$ is defined over $\Q$, $E$ has not points of order $n$ and one is interested to compute if there are points of order $n$ over a cubic field, then a necessary condition is that $\Psi_n(x)$ has an irreducible factor of degree $3$.\\
Let $E:y^2=x^3+Ax+B$ be an an elliptic curve defined over $\Q$ and $\Psi_n(x)$ its $n$-division polynomial. To determine if there exist an square free integer $d$ such that the $d$-quadratic twist of $E$ has a point of order $n$ defined over some number field $K$ it is enough to check if one the roots of $\Psi_n(x)$, say $\alpha$, belongs to $K$ and $\alpha^3+A\alpha+B=d\beta^2$ for $\beta\in K$. \\
At the Appendix appears the necessary information related to elliptic curves defined over $\Q$ with CM that it will be used to proof Theorem \ref{teo1}, \ref{teo2}, and \ref{teo3}.
\subsection{Proof of Theorem \ref{teo1}}
There are examples in Table \ref{d3} for all the cases in $\Phi^{\cm}(3)$, therefore all those torsion subgroups appear in $\Phi_{\Q}^{\cm}(3)$. This proves Theorem \ref{teo1}.
\subsection{Proof of Theorem \ref{teo2}}
It has been characterized the set $\Phi_\Q(3,G)$ for any $G\in \Phi(1)$ (see Theorem 1.2 in \cite{GJNT16}). In particular we have $\Phi^{\cm}_\Q(3,G)\subseteq \Phi_\Q(3,G)\cap \Phi_\Q^{\cm}(3)$ for any $G\in \Phi^{\cm}(1)$. Actually, except for $G=\cC_1$, the above relation is an equality since there are examples of any case in Table \ref{d3}. For trivial torsion we have $\Phi_\Q(3,\cC_1)\cap \Phi_\Q^{\cm}(3)=\{\cC_1,\cC_2,\cC_3,\cC_4,\cC_6, \cC_2\times\cC_2\}$. In Table \ref{d3} we have examples of elliptic curves $E$ with trivial torsion that over cubic fields it grows to $\cC_2,\cC_3$, and $\cC_6$. Then it remains to discard the cases $\cC_4$ and $\cC_2\times\cC_2$. In the Table \ref{isoCM} we check that if $E$ is an elliptic curve defined over $\Q$ with CM then $\mathfrak{cm}\in \{27,11,19,43,67,163\}$ or $\fcm=3$ with $E:y^2=x^3+k$ with $k\ne r^2,r^3,-432$. We split the proof depending in the cases above.
\begin{itemize}
\item $\mathfrak{cm}\in \{27,11,19,43,67,163\}$: Note that for these curves the corresponding $j$-invariants are neither $0$ nor $1728$. Then we have just quadratic twists, in particular it is only necessary to study the $n$-division polynomials for $E_\mathfrak{cm}$. In the following cases the $n$-division polynomial $\Psi_n(x)$ refers to the elliptic curve $E_\mathfrak{cm}$. We have that the field of definition of the full $2$-torsion, $\Q(E[2])$, is the splitting field of $\Psi_2(x)=f_{\mathfrak{cm}}(x)$. We have that those polynomials are irreducible and the cubic fields that they define are not a Galois extension. This proves that torsion $\cC_2\times\cC_2$ is not possible over a cubic field for those cases. In the other hand $\Psi_4(x)$ is irreducible of degree $6$ then there are not points of order $4$ over cubic field for any of the treated cases.
\item $E:y^2=x^3+k$ with $k\ne r^2,r^3,-432$: Here $\Psi_2(x)=x^3+k$ is irreducible since $k\ne r^3$, and the cubic field that it defines never is a Galois extension for any $k$. Now $\Psi_4(x)=2(x^6+20kx^3-8k^2)$, and $z=-(10\pm 6\sqrt{3})k$ is a root of $\Psi_4(\sqrt[3]{x})$. But $z=x^3$ never occurs for $x$ in a cubic field. We have proved that there are neither points of order $4$ nor full $2$-torsion over cubic fields.
\end{itemize}
This finishes the first part of the proof of Theorem \ref{teo2}. The second part is a direct consequence of the classification obtained above. We have examples at Table \ref{d3} for any set in $\mathcal{H}_\Q(3,G)$ such that all its elements belong to $\Phi^{\cm}_\Q(3,G)$. This completes the proof of Theorem \ref{teo2}. 

\subsection{Proof of Theorem \ref{teo3}}
We are going to prove Table \ref{d3}. Let $E$ be an elliptic curve defined over $\Q$ with CM. We have an explicit description at Table \ref{isoCM} of $E(\Q)_{\tors}$ in terms of its CM-invariants. Now thanks to the classification of $\Phi^{\cm}_\Q(3,G)$ for any $G\in \Phi^{\cm}(1)$ we know the possible torsion growth over cubic fields. In this case we only need to compute the $n$-division polynomials for $n\in\{2,3,7,9\}$ and check if they have (irreducible) factors of degree $3$.

First note that the torsion growth over a cubic field can only be cyclic by Theorem \ref{teo2}. Moreover, if the torsion over $\Q$ has odd order, then the $2$-division polynomial $\Psi_2(x)$ is irreducible of order $3$. Let $\alpha$ be a root of $\Psi_2(x)$ and define $K=\Q(\alpha)$. Then over $K$ the torsion is cyclic of even order. \\
We split the proof depending if the twists are quadratic or not. That is, depending if $\mathfrak{cm}\notin\{3,4\}$ or not. Suppose $\mathfrak{cm}\notin\{3,4\}$ and let $\Psi_n(x)$ denotes the $n$-division polynomial of $E_\mathfrak{cm}$.\\

\indent $\bullet$ $\fcm\in \{11,19,43,67,163\}$. The torsion over $\Q$ is trivial, therefore the torsion can grow to $\cC_2,\cC_3$ or $\cC_6$. We have that all the irreducible factor of $\Psi_3(x)$ are of even order, then no points of order $3$ over cubic fields. Only torsion growth to $\cC_2$ over the cubic field $\Q(\alpha)$, where $\Psi_2(\alpha)=0$.\\
\indent $\bullet$ $\fcm=8$. We have $E_8^k(\Q)_{\tors}\simeq \cC_2$ and $\Phi^{\cm}_\Q(3,\cC_2)=\{\cC_2,\cC_6,\cC_{14}\}$. Therefore we only need to check if $\Psi_3(x)$ and $\Psi_7(x)$ have irreducible factors of degree $3$. Again all the factors are of even degree. Then no torsion growth over cubic fields.\\
\indent $\bullet$ $\fcm\in\{ 7,28\}$. Again $E_\fcm^k(\Q)_{\tors}\simeq \cC_2$. In both cases $\Psi_3(x)$ is irreducible (of degree $4$), then no points of order $3$ over cubic fields; and $\Psi_7(x)$ has only a degree $3$ factor. In particular, these factors define cubic fields $\Q(\beta)$ that are isomorphic to $\Q(\alpha)$, where $\alpha^3+\alpha^2-2\alpha-1=0$. 
\begin{itemize}
\item For $\fcm= 7$: $\beta=36\alpha - 9$ and $f_7(\beta)=-7 (2^2 3^3 \alpha)^2$. That is, only for $k=-7$ we have points of order $7$ over a cubic field. 
\item For $\fcm =28$:  $\beta=4\alpha^2 - 4\alpha + 13$ and $f_{28}(\beta)=7(4(-3\alpha^2 + 3\alpha + 1))^2$. In this case only for $k=7$.
\end{itemize}
\indent $\bullet$  $\mathfrak{cm}=16$: For $k=1,2$ we have not torsion growth over a cubic field since for those values $E_{16}^k(\Q)_{\tors}\simeq \cC_4$. Now suppose $k\ne 1,2$. Then $E_{16}^k(\Q)_{\tors}\simeq \cC_2$. We have that there is not torsion growth over cubics since $\Psi_3(x)$ and $\Psi_7(x)$ are irreducible of degrees $4$ and $24$ respectively.\\
\indent $\bullet$ $\mathfrak{cm}=27$: Let $k=1$, then $E^1_{27}(\Q)_{\tors}\simeq \cC_3$ and $\Phi^{\cm}_\Q(3,\cC_3)=\{\cC_3,\cC_6,\cC_9\}$. We have that the torsion growth to $\cC_6$ and $\cC_9$ over $\Q(\sqrt[3]{2})$ and $\Q(\alpha)$, where $\alpha^3-3\alpha-1=0$, respectively. Now suppose $k\ne 1$ then $E^k_{27}(\Q)_{\tors}\simeq \cC_1$. There is a degree $3$ irreducible factor of $\Psi_3(x)$ such that if $\alpha$ is a root of this factor, then $\alpha=-4(2\sqrt[3]{9} + 3\sqrt[3]{3} + 1)$. Since $f_{27}(\alpha)=-3(4(4\sqrt[3]{9}+6\sqrt[3]{3}+9))^2$ we have that there are points of order $3$ over a cubic field if and only if $k=-3$ and the cubic field is $\Q(\sqrt[3]{3})$. In the other hand, the torsion growth to $\cC_2$ over $\Q(\sqrt[3]{2})$ for any $k$. \\
\indent $\bullet$ $\mathfrak{cm}=12$: For $k=1$ we have not torsion growth over a cubic field since $E_{12}^1(\Q)_{\tors}\simeq \cC_6$. Let $k\ne 1$, then $E^k_{12}(\Q)_{\tors}\simeq \cC_2$. There are not torsion growth over a cubic field to $\cC_{14}$ since all the irreducible factor of $\Psi_7(x)$ are of degree divisible by $6$.  Now the $3$-division polynomial $\Psi_3(x)$ satisfies $\Psi_3(\alpha)=0$ where $\alpha=-2\sqrt[3]{4} - 2\sqrt[3]{2} - 1$. In this case we have $f_{12}(\alpha)=-3(2(\sqrt[3]{4}+\sqrt[3]{3}+1))^2 $. That is, there are points of order $3$ over a cubic field $K$ if and only if $k=-3$ and $K=\Q(\sqrt[3]{2})$.

\

Finally the non-quadratic twists:

\

\indent $\bullet$ $\mathfrak{cm}=4$. For $k=4$ and $k=-r^2$ the torsion subgroup over $\Q$ is isomorphic to $\cC_4$ and $\cC_2\times\cC_2$ respectively. Therefore for those values there are not torsion growth over cubic fields. Suppose $k\ne4,-r^2$, then $E_4^{k}(\Q)_{\tors}\simeq \cC_2$. Then the torsion can grow over a cubic field to $\cC_6$ or $\cC_{14}$. Let $\Psi_3(x)$ and $\Psi_7(x)$ the $3$- and $7$-division polynomial, respectively, of  $E_4^{k}$. Then:
\begin{itemize}
\item $\Psi_3(x)=k^2 f_3(x^2/k)$, where $f_3(x)=3 x^2+6 x-1$ is irreducible. 
\item $\Psi_7(x)=k^{12} f_7(x^2/k)$, where $f_7(x)=7 x^{12} + 308 x^{11} - 2954 x^{10} - 19852 x^9 - 35231 x^8 - 82264 x^7 - 111916 x^6 - 42168 x^5 + 15673 x^4 + 14756 x^3 + 1302 x^2 + 196 x - 1$ is irreducible. 
\end{itemize} 
Then there can not be points of order $3$ or $7$ over cubic fields. We have proved that for the family of curves with $\mathfrak{cm}=4$ there is not torsion growth over cubic fields. \\

\indent $\bullet$ $\mathfrak{cm}=3$. In this case the elliptic curve is called Mordell curve and has the model $E_3^k:y^2=x^3+k$ for $k\in \Q^*/(\Q^*)^{6}$. Note that this case has been studied by Dey and Roy \cite{DR19}, although they used different techniques. We split the proof depending on the torsion over $\Q$:
\begin{itemize}
\item $E_3^{k}(\Q)_{\tors}\simeq \cC_6$, then $k=1$ and there are not torsion growth over cubic fields.
\item $E_3^{k}(\Q)_{\tors}\simeq \cC_3$, then $k=-432$ or $k=r^2\ne 1$. Here the torsion grows to $\cC_6$ over $\Q(\sqrt[3]{k})$, since the $2$-division polynomial is $x^3+k$ and $k$ is not a cube in $\Q$. The other possible torsion growth over a cubic is $\cC_9$. First let $k=-432$, then $g(x)=x^3 + 36 x^2 - 1728$ is the unique degree $3$ irreducible factor of the $9$-division polynomial of $E_{3}^{-432}$. Let $\alpha$ be a root of $g(x)$, then $\alpha^3-432$ is not an square in $\Q(\alpha)$. Then there is not torsion growth over $\Q(\alpha)$. Now suppose $k=r^2\ne 1$ and $P_3=(0,r)$ a point of order $3$ over $\Q$. Then $P_9=(\beta,r\gamma)\in \Q(\alpha,\beta)$ satisfies $3P_9=P_3$, where $\alpha^3-3\alpha-1=0$, $\gamma=2\alpha^2 - 4\alpha- 1$, and $\beta^3-r^2\gamma^2+r^2=0$. Therefore, in principle, the field of definition of $P_9$ is of degree $9$. We are going to check in which conditions this field is of degree $3$. Equivalently, when there is torsion growth to $\cC_9$ over a cubic field.  We need that $\beta\in \Q(\alpha)$. Note that $\beta^3=r^2(\gamma^2-1)=4(\alpha^2 - \alpha - 1)^3 r^2$. In other words, the equation $z^3=4r^2$ has solutions over $\Q( \alpha)$. But this only happens if and only if $r=4s^3$, $s\in \Q$; and $k=16$ is the unique possibility since $k$ must belong to $\Q^*/(\Q^*)^{6}$. 
\item $E_3^{k}(\Q)_{\tors}\simeq \cC_2$, then $k=r^3\ne 1$. In this case $E_3^{k}$ is the $r$-quadratic twist of $E_3$. Let $\Psi_n(x)$ be the $n$-division polynomial of $E_3$. In this case the torsion can grow over a cubic field to $\cC_6$ or $\cC_{14}$. The last case is not possible since all the irreducible factor of $\Psi_7(x)$ are of degree divisible by $6$. In the other hand $\Psi_3(x)=3x(x^3+4)$ and $f_3(\sqrt[3]{4})=-3$. Then, there are points of order $3$ over a cubic field $K$ if and only if $r=-3$ (i.e. $k=-27$) and $K=\Q(\sqrt[3]{2})$.
\item $E_3^{k}(\Q)_{\tors}\simeq \cC_1$, then $k\ne r^2,r^3,-432$. We have $\Phi_\Q^{\cm}(3,\cC_1)=\{\cC_1,\cC_2,\cC_3,\cC_6\}$. We are going to study the $n$-division polynomial, $\Psi_n(x)$, of $E_3^{k}$:
\begin{itemize}
\item $\Psi_2(x)=x^3+k$ is irreducible, then there is a point of order $2$ over $\Q(\sqrt[3]{k})$.
\item $\Psi_3(x)=3x(x^3+4k)$. Note that if $x=0$ then the equation $y^2=k$ has solution over a cubic field if and only if $k$ is an square over $\Q$. But we have assumed that $k\ne r^2$. Let $\alpha\ne 0$ be another root of $\Psi_3(x)=0$. Then $y^2=\alpha^3+k=\alpha^3+4k-3k=-3k$ has solution over a cubic field if and only if $k=-3s^2$ for some $r\in \Q$. In particular the cubic field is $\Q(\sqrt[3]{12s^2})$.
\end{itemize}
Finally we study the torsion growth over a cubic field $K$ to $\cC_6$. Necessary $k=-3s^2$ and the cubic fields of definition of the points of order $2$ and $3$ must be equal to $K$. From the equality $\Q(\sqrt[3]{3s^2})=\Q(\sqrt[3]{12s^2})$ we obtain $K=\Q(\sqrt[3]{4})$. In the other hand, $\sqrt[3]{3s^2}\in K$ if and only if $s=6t^3$; but necessarily $t=\pm 1$ since $k\in \Q^*/(\Q^*)^{6}$. Then we finish that the torsion growth over a cubic field $K$ to $\cC_6$ if and only if $k=-108$ and $K=\Q(\sqrt[3]{2})$.
\end{itemize}

{\bf Remark:} {\em All the computation have been done using \texttt{Magma} \cite{magma} and the source code is available in the online supplement \cite{MagmaCode}.}

{
\section*{Appendix. Elliptic curve over $\Q$ with CM.}\label{appendix}
The necessary information related to elliptic curves with CM to be used in this paper appear in this Appendix. Let $E$ be an elliptic curve defined over $\Q$ with CM  by an order $R=\Z+\mathfrak{f}\,\mathcal{O}_K$ of conductor $\mathfrak{f}$ in a quadratic imaginary field $K=\Q(\sqrt{-D})$, where $ \mathcal{O}_K$ is the ring of integer of $K$. Then $R$ is one of the thirteen orders that correspond to the first and second column of Table \ref{isoCM}. Each order correspond to a $\Qbar$-isomorphic class of elliptic curves defined over $\Q$ with CM. The corresponding $j$-invariant appears at the third column. Fourth column, $\mathfrak{cm}$, denotes the absolute value of the discriminant of the CM quadratic order $R$. Note that the integer $\fcm$ gives the $\Qbar$-isomorphic class of $E$. Fifth column gives a pair of integers $[A_\fcm,B_\fcm]$ such that if we denote by $f_\fcm(x)=x^3+A_\fcm x+B_\fcm$ then $E_\fcm:y^2=f_\fcm(x)$ is an elliptic curve with $j(E_\fcm)$ equal to the $j$-invariant $j$ at the same row.  That is, $E_\fcm$ is a representative for each class. Now by the theory of twists of elliptic curves (cf. \cite[X \S 5]{Silverman}) applied to elliptic curve defined over $\Q$ with CM we have:
\begin{itemize}
\item If $\fcm\in\{12,27,16,7,28,11,19,43,67,163\}$ (i.e. $j(E)\ne 0,1728$) then $E$ is $\Q$-isomorphic to the $k$-quadratic twist of $E_\fcm$ for some squarefree integer $k$. That is, $E$ has a short Weierstrass model of the form $E_\fcm^k:y^2=x^3+k^2A_\fcm x+k^3B_\fcm$.
\item If $\fcm=3$ (i.e. $j(E)=0$) then $E$ has a short Weierstrass model of the form $E_3^k:y^2=x^3+k$, where $k$ is an integer such that $k\in\Q^*/(\Q^*)^{6}$.
\item If $\fcm=4$ (i.e. $j(E)=1728$) then $E$ has a short Weierstrass model of the form $E_4^k:y^2=x^3+kx$, where $k$ is an integer such that $k\in\Q^*/(\Q^*)^{4}$.
\end{itemize}
Note that $k$ and $\fcm$ are uniquely determined by $E$. We call them the CM-invariants of the elliptic curve $E$.

\

Finally, given an elliptic curve $E$ defined over $\Q$ with CM, at the last two columns of Table \ref{isoCM} we give a characterization of its torsion subgroup (over $\Q$) depending on its CM-invariants $(\fcm,k)$ (see Table 3 at \cite[\S 2]{GJCM2}).

\begin{table}[ht!]
\caption{Elliptic curves defined over $\Q$ with CM. Torsion over $\Q$.}\label{isoCM}
\begin{tabular}{|c|c|c|c|c|c|c|}
\hline
$-D$ & $\mathfrak{f}$ &  $j$ & $\mathfrak{cm}$ & $[A_{\mathfrak{cm}},B_{\mathfrak{cm}}]$ & $k$ & $E^k_{\mathfrak{cm}}(\Q)_{\tors}$  \\
\hline
\multirow{8}{*}{$-3$}  & \multirow{4}{*}{$1$}  & \multirow{4}{*}{$0$}  & \multirow{4}{*}{$3$}  & \multirow{4}{*}{$$[0,1]$$}  & $1$ & $\cC_6$ \\
\cline{6-7}
& & & & &  $-432,r^2\ne 1$&  $\cC_3$\\
\cline{6-7}
& & & & &  $r^3 \ne 1$&  $\cC_2$\\
\cline{6-7}
& & & & &  $\ne r^2,r^3,-432$ &  $\cC_1$ \\
\cline{2-7}
 & \multirow{2}{*}{$2$}  &  \multirow{2}{*}{$2^4\cdot 3^3\cdot 5^3$} & \multirow{2}{*}{$12$} & \multirow{2}{*}{$[-15,22]$} & $1$ & $\cC_6$ \\
\cline{6-7}
& & & & &   $\ne 1$ & $\cC_2$ \\
\cline{2-7}
& \multirow{2}{*}{$3$}  &  \multirow{2}{*}{$-2^{15}\cdot  3\cdot  5^3$} & \multirow{2}{*}{$27$} & \multirow{2}{*}{$[-480,4048]$} & $1$ & $\cC_3$ \\
\cline{6-7}
& & & & &   $\ne 1$ & $\cC_1$ \\
\hline
\multirow{5}{*}{$-4$} & \multirow{3}{*}{$1$} & \multirow{3}{*}{$2^6\cdot  3^3=1728$} & \multirow{3}{*}{$4 $} & \multirow{3}{*}{$[1,0]$} & $4$ & $\cC_4$ \\
\cline{6-7}
& & & & &   $-r^2$ & $\cC_2\times \cC_2$ \\
\cline{6-7}
& & & & &   $\ne 4, -r^2$ & $\cC_2$ \\
\cline{2-7}

& \multirow{2}{*}{$2$ }& \multirow{2}{*}{$2^3\cdot  3^3\cdot  11^3$} & \multirow{2}{*}{$16 $} &  \multirow{2}{*}{$[-11,14]$} & $1,2$ & $\cC_4$ \\
\cline{6-7}
& & & & &   $\ne 1,2$ & $\cC_2$ \\

\hline
 \multirow{2}{*}{$-7$} & $1$ & $-3^3\cdot  5^3$ & $7 $ &$[- 2835,- 71442]$ & $-$ & $\cC_2$ \\
\cline{2-7}
 & $2$ & $3^3\cdot  5^3\cdot  17^3$ & $28 $ &$  \quad [-595,5586]$ & $-$ & $\cC_2$ \\
\hline
 $-8$ & $1$ & $2^6\cdot  5^3$ & $ 8$ &$[- 4320, 96768]$ & $-$ & $\cC_2$ \\
\hline
 $-11$ & $1$ & $-2^{15}$ &$ 11$ & $ [ - 9504 , 365904]$ & $-$ & $\cC_1$ \\
\hline
 $-19$ & $1$  & $-2^{15}\cdot  3^3$ & $19 $ & $ [ - 608 , 5776] $ & $-$ & $\cC_1$ \\
\hline
 $-43$ & $1$ & $-2^{18}\cdot  3^3\cdot  5^3$& $43 $ & $  [ - 13760, 621264]$  & $-$ & $\cC_1$ \\
\hline
 $-67$ & $1$  & $2^{15}\cdot  3^3\cdot  5^3\cdot  11^3$& $67 $ & $  [- 117920 , 15585808] $ & $-$ & $\cC_1$ \\
\hline
 $-163$ & $1$  & $-2^{18}\cdot  3^3\cdot  5^3\cdot  23^3\cdot  29^3$ & $163 $ & $  [- 34790720 , 78984748304]  $ & $-$ & $\cC_1$ \\
\hline
\end{tabular}

\end{table}

}

\end{document}